\documentstyle[12pt]{article}
\advance\textheight by 0.5cm
\advance\voffset by -0.2cm
\makeatletter
\ifnum\@ptsize=2  \fi
\ifnum\@ptsize=0  \fi
\if@twoside\fi
\def\runningtitle#1{%
   \global\def\rightmark{%
      \rlap{\makebox[\textwidth]{\centering\small #1}}%
   }
}
\let\oldtitle=\title
\let\oldauthor=\author
\let\olddate=\date
\let\oldmaketitle=\maketitle
\def\title#1{\def\thetitle{#1}\oldtitle{\vskip-1cm\Large\bf #1}}
\def\date#1{\def\thedate{#1}\olddate{\normalsize #1}}
\def\author#1{\def\theauthor{#1}\oldauthor{\normalsize #1}}
\date{\today} 
\def\maketitle{\oldmaketitle\thispagestyle{empty}}
\def\tighttoc{
   \def\l@section{\@dottedtocline{1}{0em}{1.6em}}
   \def\l@subsection{\@dottedtocline{2}{1.5em}{3.2em}}
}
\tighttoc
\def\tightbib{
   \let\oldthebibliography=\thebibliography
   \def\thebibliography##1{
      \small
      \def\@listi{\topsep=0cm \parsep=0cm \itemsep=0cm}
      \oldthebibliography{##1}
   }
}
\tightbib
\def\@sect#1#2#3#4#5#6[#7]#8{\ifnum #2>\c@secnumdepth
     \let\@svsec\@empty\else
     \refstepcounter{#1}\edef\@svsec{\csname the#1\endcsname. }\fi
     \@tempskipa #5\relax
      \ifdim \@tempskipa>\z@
        \begingroup #6\relax
          \@hangfrom{\hskip #3\relax\@svsec}{\interlinepenalty \@M #8\par}%
        \endgroup
       \csname #1mark\endcsname{#7}\addcontentsline
         {toc}{#1}{\ifnum #2>\c@secnumdepth \else
                      \protect\numberline{\csname the#1\endcsname.}\fi
                    #7}\else
        \def\@svsechd{#6\hskip #3\relax  
                   \@svsec #8\csname #1mark\endcsname
                      {#7}\addcontentsline
                           {toc}{#1}{\ifnum #2>\c@secnumdepth \else
                             \protect\numberline{\csname the#1\endcsname.}\fi
                       #7}}\fi
     \@xsect{#5}}
\def\centersections{
   \def\section{\@startsection{section}{1}{\z@}%
      {2\bigskipamount plus \medskipamount minus \medskipamount}%
      {2.3ex plus.2ex}{\centering\reset@font\normalsize\bf\boldmath}}
   \def\subsection{\@startsection{subsection}{2}{\z@}%
      {\bigskipamount plus \medskipamount minus \smallskipamount}%
      {1.5ex plus.2ex}{\centering\reset@font\normalsize\it}}
}
\centersections
\newlength{\paritemlabelwidth}
\def\paritemlabel{$\bullet$}
\def\paritems#1{
   \def\@listi{\topsep=0cm\parsep=0cm\itemsep=0cm\partopsep=\medskipamount}
   \let\@listii=\@listi
   \let\@listiii=\@listi
   \@listi
   \def\paritemarg{#1}
   \ifx\paritemarg\empty\relax\else\def\paritemlabel{\paritemarg}\fi
   \begin{trivlist}\item[]\indent\mbox{}\vskip-\baselineskip
      \settowidth{\paritemlabelwidth}{\paritemlabel}%
      \def\paritem{\@ifnextchar[{\bparitem}{\nparitem}}
      \def\bparitem[##1]{
          \par\makebox[\paritemlabelwidth][r]{##1}\hspace{1ex}\ignorespaces}
      \def\nparitem{\bparitem[\paritemlabel]}
      \let\item=\paritem
      \renewcommand{\i}{\item[i)]}
      \newcommand{\ii}{\item[ii)]}
      \newcommand{\iii}{\item[iii)]}
}
\def\endparitems{\end{trivlist}}
\def\bbC{{\mathchoice {\setbox0=\hbox{$\displaystyle\rm C$}\hbox{\hbox
   to0pt{\kern0.4\wd0\vrule height0.9\ht0\hss}\box0}}
   {\setbox0=\hbox{$\textstyle\rm C$}\hbox{\hbox
   to0pt{\kern0.4\wd0\vrule height0.9\ht0\hss}\box0}}
   {\setbox0=\hbox{$\scriptstyle\rm C$}\hbox{\hbox
   to0pt{\kern0.4\wd0\vrule height0.9\ht0\hss}\box0}}
   {\setbox0=\hbox{$\scriptscriptstyle\rm C$}\hbox{\hbox
   to0pt{\kern0.4\wd0\vrule height0.9\ht0\hss}\box0}}}}
\def\bbP{{\rm I\!P}}
\def\TheoremIndent{\hspace*{\parindent}}
\def\TheoremNumber#1{(#1) }
\def\@begintheorem#1#2{\trivlist%
   \item[\hskip\labelsep{\TheoremIndent\bf\TheoremNumber{#2}#1.}]\it%
}
\def\@opargbegintheorem#1#2#3{\trivlist
   \item[\hskip\labelsep{\TheoremIndent\bf\TheoremNumber{#2}#1 \rm(#3).}]\it%
}
\newtheorem{satz}{Satz}[section]
\newtheorem{theorem}[satz]{Theorem}

\newtheorem{lemma}[satz]{Lemma}
\newtheorem{corollary}[satz]{Corollary}
\newtheorem{proposition}[satz]{Proposition}

\newenvironment{proclaim*}[1]{
   \trivlist\item[]\it
   \TheoremIndent{\bf #1.}%
}{
   \endtrivlist
}

\newenvironment{assertion*}{%
   \trivlist\item[]\it%
   \TheoremIndent\ignorespaces
}{
   \endtrivlist
}
\def\startproof{\par\addvspace{\bigskipamount}}
\def\proof{\startproof{\it Proof. }}

\def\qed{\nopagebreak\hspace*{\fill}
   \frame{\rule[0pt]{0pt}{8pt}\rule[0pt]{8pt}{0pt}}
   \par\addvspace{\bigskipamount}
}

\def\to{\mathop{\longrightarrow}\limits}
\def\mapsto{\mathop{\mapstochar\longrightarrow}\limits}
\def\phi{\varphi}
\def\epsilon{\varepsilon}
\def\tilde{\widetilde}

\def\({\left(}
\def\){\right)}
\def\be{\begin{eqnarray*}}
\def\ee{\end{eqnarray*}}
\def\tfrac#1#2{{\textstyle\frac{#1}{#2}}}

\def\d#1#2{\frac{\partial #1}{\partial #2}}

\def\with{\mid}

\def\inverse{^{-1}}
\def\restr#1{\big|_{#1}}
\def\isom{\cong}

\def\rank{\mathop{\rm rank}\nolimits}

\def\PGL{\mathop{\rm PGL}\nolimits}

\def\and{\quad\mbox{ and }\quad}
\def\eqnref#1{(\ref{#1})}

\def\tfrac#1#2{{\textstyle\frac{#1}{#2}}}
\def\vect#1{\left(\begin{array}{c} #1 \end{array}\right)}
\def\eps{\epsilon}

\begin{document}

\title{Poncelet theorems}
\runningtitle{Poncelet theorems}
\author{W.\ Barth, Th.\ Bauer}
\date{ }
\maketitle
\tableofcontents


\setcounter{section}{-1}
\section{Introduction}

   The aim of this note is to collect some more or less classical
   theorems of Poncelet type and to provide them with short modern
   proofs. Where classical geometers used elliptic functions
   (or angular functions), we use elliptic curves (or degenerate elliptic
   curves decomposing into two rational curves). In this way we
   unify the geometry  underlying these Poncelet type statements.

   Our starting point is a space Poncelet theorem for
   two quadrics in $\bbP_3$ (Sect.\ \ref{section Weyr}).
   This seems to have been observed
   first by Weyr \cite{Weyr} p.\ 28, and amplified by Griffiths and
   Harris \cite{GriHar77}.

   The classical Poncelet theorem (cf.\ \cite{GriHar78a}),
   a statement on two conics in the plane $\bbP_2$,
   follows from Weyr's space Poncelet theorem, if one of
   the quadrics is taken to be a cone, see Sect.\ \ref{section Weyr}.

   Gerbaldi \cite{Ger19} gave formulas counting the
   number of conics in a pencil which are in Poncelet position
   with respect to a fixed conic in the pencil. We show that
   his formulas are simple consequences of the space
   Poncelet theorem (Sect.\ \ref{section Gerbaldi}).

   In Sect.\ \ref{section revolution}
   we evaluate explicitly the space Poncelet
   condition for two quadrics of revolution about the
   same axis.

   We show that theorems such as
   Emch's theorem on circular series \cite{Emc01}
   and (a complex-projective version of) the
   'zig-zag' theorem \cite{BlaHowHow74}
   can be understood
   by considering torsion points on
   elliptic curves
   (see Sects.\ \ref{section Emch Steiner} and
   \ref{section zig-zag}).
   Further, we prove a Poncelet version of the Money-Coutts theorem
   (Sect.\ \ref{section Money-Coutts}).

   Although for most of the contents of this paper only the
   presentation is new, it seems worth--while to us to consider
   Poncelet type theorems from a modern geometric point of view.
   In this spirit the equations of modular curves given in
   \cite{BarMic93} have been transformed
   by N. Hitchin \cite{Hitchin} into solutions of
   the Painlev\'{e}--VI--equation. And E.\ Previato relates Poncelet
   theorems to integrable Hamiltonian systems and billiards
   \cite{Prev}.

\bigskip
   {\em Conventions.} The base field always is the field
   $\bbC$ of complex numbers.

   If we mention circles, quadrics of revolution, or spheres, we
   mean the corresponding varieties over $\bbC$.

\bigskip
   {\em Acknowledgement. This research was supported by DFG contract
   Ba 423/7-1 and EG contract SC1--0398--C(A).
   }


\section{Weyr's Poncelet theorem in $\bbP_3$}
\label{section Weyr}

   Let $Q_1,Q_2\subset\bbP_3$ be quadrics of ranks $\ge 3$ such that
   their intersection curve $E=Q_1 \cap Q_2$ is either a smooth
   elliptic curve or the union of two conics $C_1,C_2$ meeting in
   two distinct points.  We fix rulings $R_1$ on $Q_1$ and $R_2$ on $Q_2$.

\begin{theorem}[\cite{Weyr} p.\ 28, \cite{Hurw} p.\ 13]
   \label{3-dim Poncelet}
   Suppose that there exists a closed sequence of distinct lines
   $L_1,\dots,L_{2n},L_{2n+1}=L_1$ such that the line
   $L_i$ belongs to $R_1$ resp.\ $R_2$, if $i$ is odd resp.\ even,
   and such that consecutive lines $L_i,L_{i+1}$ intersect each other.
   Then there are such closed sequences of length $2n$ through
   any point on $Q_1\cap Q_2$.
\end{theorem}

\proof
   The rulings $R_1,R_2$ define involutions $\iota_1,\iota_2$
   on $E$
   interchanging the two intersection points of $E$ with a line in
   $R_1$ resp.\ $R_2$.  Let $t:E\to E$ be the composition $\iota_2\iota_1$
   and let $e:=L_{2n}\cap L_1$.
   Then $L_{2k-1}$ is the unique line in $R_1$ through
   $t^{k-1}e$ and
   $L_{2k}$ is the unique line in $R_2$ through
   $\iota_1t^{k-1}e$ for $k\ge 1$.
   The closedness $L_{2n+1}=L_1$ implies $t^n(e)=e$.

   First we consider the case that $E$ is smooth elliptic.
   The involutions $\iota_1,\iota_2$ have fixed points, so their
   composition $t$ is a translation on $E$.
   But then $t^n(e)=e$ and the assumption that the lines $L_1,\dots,L_n$
   are distinct
   imply that $t$ is of order $n$, hence the
   assertion.

   Now assume $E=C_1\cup C_2$ is the union of two conics.
   The involutions $\iota_1$ and $\iota_2$ interchange these
   conics. So $t(C_i)=C_i, \, i=1,2$ and the two intersections
   of $C_1$ and $C_2$ are fixed points of
   $t$.
   Therefore $t$ induces two automorphisms
   $
    t_i:\bbC^*\isom C_i\setminus(C_1\cap C_2)
    \to              C_i\setminus(C_1\cap C_2)\isom\bbC^*
   $.
   If $e$ lies on $C_1$, say, then
   $t^n(e)=e$ implies that $t_1$ is the multiplication by a primitive
   $n$-th root of unity, i.e.\ its order is $n$.
   But then $t_2$ is of order $n$ as well, because of
   $t_2=\iota_2 t_1\inverse\iota_2$.
\qed

   Theorem \ref{3-dim Poncelet}
   also holds for two different rulings $R_1$ and $R_2$
   on the same smooth quadric $\simeq \bbP_1 \times \bbP_1$, if
   $E$ is a curve on this quadric of bidegree $(2,2)$, either
   smooth elliptic or the union of two rational curves meeting
   in two distinct points.

   The proof is literally the same.

   We say that the pair of rulings $(R_1,R_2)$ {\it satisfies the
   Poncelet-$n$-condition} if the automorphism $t:E\to E$ defined above
   is of order $n$.  A priori this is a property of the {\it ordered}
   pair $(R_1,R_2)$.  But the automorphism associated to $(R_2,R_1)$ is
   just $t^{-1}$, thus:

\begin{assertion}
   The pair of rulings $(R_1,R_2)$ satisfies the Poncelet-$n$-condition
   if and only if $(R_2,R_1)$ does.
\end{assertion}

   Already Hurwitz \cite{Hurw}, p.\ 13, observed that Theorem
   \ref{3-dim Poncelet}
   implies the usual Poncelet theorem in $\bbP_2$.

   In fact, if one of the quadrics,
   $Q_1$, is smooth and the other one, $Q_2$, is a cone
   with its top $P_0$ not on $E=Q_1 \cap Q_2$ and $E$ is
   smooth elliptic, then
   the three-dimensional Poncelet theorem is equivalent
   to Poncelet's theorem for two conics in the plane.  To see this,
   we denote by $\pi_i:Q_i\to\bbP_2$ the projections
   from $P_0$.  The morphism $\pi_1$ is of degree 2, ramified over a
   smooth conic $C\subset\bbP_2$.  The image of $\pi_2$ is another
   smooth conic $D\subset\bbP_2$ in general position with respect
   to $C$.  We have:

\begin{proposition}\label{23}
   The quadrics $Q_1$ and $Q_2$ are in Poncelet-$n$-position if and
   only if the conic $C$ is $n$-inscribed into $D$.
\end{proposition}

   We say: $C$ is $n$--inscribed into $D$, if there is a polygon
   consisting of $n$ tangents to $C$ with its vertices on
   $D$, \cite{BarMic93}.

\proof
   We choose an isomorphism $Q_1\isom\bbP_1\times\bbP_1$ such that
   the ruling $R_1$ is parametrized by the first factor.
   Let $\Delta = \pi^{-1}C \subset \bbP_1 \times \bbP_1$ be the
   ramification curve. It is a curve of bidegree $(1,1)$
   on $\bbP_1 \times \bbP_1$ and
   induces an isomorphism between both copies of $\bbP_1$
   such that $\Delta$ is the diagonal.
   We consider the
   map
   \be
      Q_1   &\to&   \bbP_2\times C^*   \\
      (a,b) &\mapsto&  (u,T)
   \ee
   where $u:=\pi_1(a,b)$ and $T:=\pi_1(a\times\bbP_1)$.  So $T$ is the
   tangent from $u$ to $C$ meeting $C$ in the point $\pi_1(a,a)$.
   The projection $\pi_1$ restricts to a map $E\to D$ of degree 2.
   Over a point $u\in D$ we have points $(a,b),(b,a)\in E$.
   Now we determine the map on the pairs $(u,T)$ induced by $t$.
   We start with a smooth point $P=(a,b)$ on $E$.
   The line $L_1=a\times\bbP_1\in R_1$ through
   $P$ meets $E$ in another point $P'=(a,b')$.  This means that we
   pass from the pair $(u,T)$ to $(u',T)$, where $u':=\pi_1(a,b')$ is
   the second intersection point of $T$ and $D$.
   The image point $P'':= t(P)$ is the second intersection point of the
   line $L_2\in R_2$ through $P'$ with $E$.  Thus we have $P''=(b',a)$,
   because points in the same ray of $Q_2$ lie in the same fibre of the
   projection $\pi_2$.  Passing from the point $(a,b')$ to $(b',a)$
   means passing from the pair $(u',T)$ to $(u',T')$, where $T'$ is
   the second tangent to $C$ through the point $u'$.
   So we see that the map
   $$
      t:(u,T)\mapsto (u',T')
   $$
   just describes the Poncelet process for the pair of conics $C,D$.
   This implies the assertion.
\qed

   We should point out here, how Theorem \ref{3-dim Poncelet}
   relates to the Poncelet theorem in space considered by
   Griffiths and Harris in \cite{GriHar77} on finite polyhedra
   both inscribed and circumscribing two smooth quadrics $Q_1$ and $Q_2$:
   Such a polyhedron exists if and only if both of the pairs
   $R_1,R_2$ and $R_1',R_2'$ satisfy a Poncelet condition in the sense
   above, where $R_i,R_i'$ are the rulings on the dual quadrics
   $Q_i^*$, $i=1,2$.


\section{Poncelet pairs in a pencil}

   Let $Q_1,Q_2$ be quadrics of rank $\ge 3$ such that
   the pencil $\lambda_1 Q_1+\lambda_2 Q_2$ generated by them is generic.
   This means that its base locus $E = Q_1\cap Q_2$
   is smooth, or equivalently that
   the discriminant $d(\lambda_1,\lambda_2)=det(\lambda_1 Q_1+\lambda_2 Q_2)$
of
   the pencil has no multiple roots.
   The four roots of $d(\lambda_1,\lambda_2)$
   correspond to the cones in the pencil,
   i.e.\ to those quadrics carrying only one ruling.
   Let $M\to\bbP_1$ be the double cover of $\bbP_1$ branched over the
   roots of $d(\lambda_1,\lambda_2)$.  The points of the elliptic curve $M$
   can be identified with the rulings on the quadrics of the given pencil.
   Each ruling $R$ in $M$ defines an involution $I_R:E\to E$ with fixed
   points.  Choosing an origin in $E$,
   we can write $I_R$ as $x\mapsto -x+a$
   with a unique point $a\in E$, so we obtain a map
   \be
      \Phi:M&\to& E \\
      R &\mapsto& a \ .
   \ee

\begin{proposition}
   $\Phi$ is an isomorphism of groups (if the origin of $M$ is chosen
   appropriately).
\end{proposition}

\proof
   It suffices to show that $\Phi$ is injective.  So let $R_1,R_2\in M$
   such that $\Phi(R_1)=\Phi(R_2)$, i.e.\ such that the involutions
   $I_{R_1},I_{R_2}$ coincide.  Let $P\in E$ be a point which is not
   fixed under $I_{R_1}$ or $I_{R_2}$.  Because of $I_{R_1}=I_{R_2}$
   the lines
   \be
      \overline{P,I_{R_1}(P)}\in R_1 \\
      \overline{P,I_{R_2}(P)}\in R_2
   \ee
   are equal.  Varying the point $P$ in $E$ we conclude $R_1=R_2$.
\qed

   This gives the following corollary which will be applied in
   Sect.\ \ref{section Gerbaldi}.

\begin{corollary}\label{torsion criterion}
   a) Two rulings $R_1,R_2$ in the pencil satisfy the Poncelet-$n$-condition,
   if and only if the point $R_2-R_1\in M$ is a primitive $n$-torsion
   point.

   b) For a fixed ruling $R_1\in M$ there are $T(n)$ rulings $R_2\in M$ such
   that $R_1,R_2$ satisfy the Poncelet-$n$-condition.
   Here $T(n)$ denotes the number of primitive $n$-torsion points on
   an elliptic curve.
\end{corollary}

   This corollary is independent of the choice of origin on $M$,
   if one interprets $R_2-R_1$ as the translation on $M$ mapping $R_1$
   to $R_2$.

   In general the Poncelet property depends on the choice of the rulings
   $R_1,R_2$ on $Q_1,Q_2$.  However
   if it holds for $R_1$ and $R_2$, then it also holds for
   the complementary rulings $R_1'$ on $Q_1$ and $R_2'$ on $Q_2$.
   In fact, $R_2-R_1=-(R_2'-R_1')$, independently of the choice of
   origin on $M$. This is even true if one of the quadrics,
   say $Q_2$ is a cone and $R_2'=R_2$:

\begin{assertion}
   Let $Q_1$ be smooth and let $Q_2$ be a cone.  We denote by
   $R_1,R'_1$ the rulings on $Q_1$ and by $R_2$ the unique ruling
   on $Q_2$.  Then $R_1',R_2$ satisfy the Poncelet-$n$-condition
   if and only if $R_1,R_2$ do.
\end{assertion}


\section{Gerbaldi's formula for the number of inscribed conics}
\label{section Gerbaldi}

   Gerbaldi \cite{Ger19} considered the invariant of two conics
   $C$ and $D \subset \bbP_2$ vanishing if $C$ is $n$--inscribed
   into $D$. Using elliptic functions and continued fractions
   he showed:

\begin{theorem}[\cite{Ger19}, p.\ 103] \label{Gerb}
   This invariant is of degree $\frac{1}{2} T(n)$ in $C$ and of degree
   $\frac{1}{4} T(n)$ in $D$.
\end{theorem}

   Here we do not want to make his assertion precise. Instead we
   prove the following two assertions, together equivalent with
   Theorem \ref{Gerb}. Ignorantly of Gerbaldi's paper they
   were shown in \cite{BarMic93} using a rational elliptic surface.

\begin{theorem}\label{number inscribed}
   Let $\lambda C+\mu D$, $(\lambda:\mu)\in\bbP_1$,
   be a generic pencil of conics in $\bbP_2$.
   Then the number of conics in the pencil, which are $n$-inscribed
   into $D$ is $\frac12 T(n)$.
\end{theorem}

\proof
   1)
   There is a smooth quadric $Q_1\subset\bbP_3$ and a cone $Q_2\subset\bbP_3$
   such that
   \begin{paritems}{iii)}
      \item[i)]
         the branch locus of the projection $Q_1\to\bbP_2$ from the top $P_0$
         of $Q_2$ is $C$,
      \item[ii)]
         the image of $Q_2$ under $P_0$-projection is $D$, and
      \item[iii)]
         the pencil generated by $Q_1$ and $Q_2$ is generic.
   \end{paritems}

   In fact, in suitable homogeneous coordinates
   $x,y,t$ on $\bbP_2$ and $x,y,z,t$ on $\bbP_3$
   the conics $C$ and $D$
   are given by equations
   $$
         C:  x^2+y^2+t^2 = 0 , \qquad
         D:  \alpha x^2+\beta y^2+\gamma t^2 = 0
   $$
   with $\alpha,\beta,\gamma\in\bbC$.
   We choose $Q_2$ to be the cone
   $$
      Q_2:  \alpha x^2+\beta y^2+\gamma t^2 = 0 \\
   $$
   with top $P_0=(0:0:1:0)$ and $Q_1$ to be the smooth quadric
   $$
      Q_1:  x^2+y^2+z^2+t^2 = 0 \ .
   $$
   Now conditions i) and ii) are obviously satisfied.  For condition iii)
   note that the discriminant of the pencil $\lambda_1 Q_1+\lambda_2 Q_2$,
   $$
      d(\lambda_1,\lambda_2)=\lambda_1(\lambda_1+\lambda_2\alpha)
          (\lambda_1+\lambda_2\beta)(\lambda_1+\lambda_2\gamma) \ ,
   $$
   has no multiple roots, since the pencil $\lambda C+\mu D$ is
   generic by assumption.

   2)
   Let $Q_1,Q_2\subset\bbP_3$ be quadrics with the properties i),ii),iii)
   above.  We claim that the branch locus of a smooth quadric
   $Q_{\lambda_1,\lambda_2}=\lambda_1 Q_1+\lambda_2 Q_2$
   under $P_0$-projection is a conic in the pencil $\lambda C+\mu D$.
   It is enough to show that the branch divisors of the
   $P_0$-projections $Q_{\lambda_1,\lambda_2}\to\bbP_2$ vary in a pencil.
   Now, a point $P\in Q_{\lambda_1,\lambda_2}$ is a branch point,
   if the line $\overline{PP_0}$ touches $Q_{\lambda_1,\lambda_2}$
   in $P$, i.e.\ if $P_0^t Q_{\lambda_1,\lambda_2} P=0$ (in matrix notation).
   So the branch divisor is just the intersection of
   $Q_{\lambda_1,\lambda_2}$ with the polar
   $P_0^t Q_{\lambda_1,\lambda_2}$ of the point $P_0$.
   But
   $$
      P_0^t Q_{\lambda_1,\lambda_2} = P_0^t(\lambda_1 Q_1+\lambda_2 Q_2)
         = \lambda_1 P_0^t Q_1
   $$
   shows that this polar is the same for all quadrics in the pencil.
   So the branch divisors vary in a pencil on this polar.

   3) If a quadric $Q_{\lambda_1,\lambda_2}$ of the pencil is smooth,
   then according to Proposition \ref{23} the quadrics
   $Q_{\lambda_1,\lambda_2}$ and $Q_2$ are in Poncelet-$n$-position,
   if and only if the branch locus of $Q_{\lambda_1,\lambda_2}$
   is $n$-inscribed into $D$.
   If $Q_{\lambda_1,\lambda_2}$ is a cone, then $Q_{\lambda_1,\lambda_2}$
   and $Q_2$ are certainly not in Poncelet-$n$-position, because
   their rulings are halfperiods on the elliptic curve $M$
   parametrizing the rulings in the pencil.
   We conclude that the number of conics in the pencil
   $\lambda C+\mu D$, which are $n$-inscribed into $D$ equals the
   number of quadrics $Q_{\lambda_1,\lambda_2}$ such that
   $Q_{\lambda_1,\lambda_2}$ and $Q_2$ are in Poncelet-$n$-position.
   But this number is just half the number of rulings $R\in M$
   such that $R,R_2$ satisfy the Poncelet-$n$-condition.
   Now the assertion of the theorem follows by
   Corollary \ref{torsion criterion}.
\qed

\begin{theorem}\label{number circumscribed}
   Let $\lambda C + \mu D$ be a general pencil of conics in $\bbP_2$.
   Then the number of conics in this pencil, $n$--circumscribed
   about $D$ is $\frac{1}{4} T(n)$.
\end{theorem}

\proof
   A conic $\lambda C+\mu D$ is $n$--circumscribed about $D$ if
   and only if the dual conic $(\lambda C + \mu D)^*$ is
   $n$--inscribed into $D^*$. The dual pencil
   $(\lambda C + \mu D)^*$ is parametrized by a smooth
   conic $\Gamma$ in the space $\bbP_5$ of all conics.

   Denote by $H \subset \bbP_5$ the hypersurface parametrizing
   conics $n$--circumscribed about $D$. By
   Theorem \ref{number inscribed}
   $$ \deg(H) = \frac{1}{2} \Gamma . H = \frac{1}{4} T(n).$$
   This shows that the general pencil $\lambda C + \mu D$ meets
   $H$ in $\frac{1}{4} T(n)$ points.
\qed


\section{A Poncelet theorem on three conics}
\label{section three conics}

   As an application of the three-dimensional Poncelet theorem
   we now prove a Poncelet theorem on three conics.

   Let $C,C_1,C_2\subset\bbP_2$ be three smooth conics in a generic pencil.
   Let $P_1$ be an arbitrary point on $C$.  There are two tangent lines
   to $C_1$ through $P_1$.  We choose one of them, $T_1$ say, and define
   $P_2$ to be its second point of intersection with $C$.
   Next, we choose a tangent $T_2$ to $C_2$ through $P_2$.
   We will now describe, how the data $P_1,T_1,T_2$ determine a sequence
   $(T_i)_{i\ge 1}$ of lines such that for $i\ge 1$ the line $T_i$ is
   tangent to $C_1$ resp.\ $C_2$, if $i$ is odd resp.\ even, and such
   that the intersection points of consecutive lines $T_i,T_{i+1}$
   lie on $C$.

   To this end, we choose
   (as in the proof of Theorem \ref{number inscribed})
   smooth quadrics $Q_1,Q_2$ and a
   cone $Q$ in a generic pencil in $\bbP_3$ such that,
   denoting the projection $\bbP_3 - - \rightarrow\bbP_2$
   from the top of $Q$ by $\pi$,
\begin{paritems}{ii)}
   \item[i)] the branch locus of $\pi\restr{Q_i}:Q_i\to\bbP_2$ is $C_i$,
   $i=1,2$, and
   \item[ii)] the image of $Q$ under $\pi$ is $C$.
\end{paritems}
   Then we have $P_1=\pi(e)$ for some point $e$ on the elliptic curve
   $E=Q_1\cap Q_2$.  The tangent $T_1$ is the image $\pi(L_1)$ of a line
   $L_1$ on $Q_1$ through $e$.  Let $R_1$ be the ruling on $Q_1$ containing
   $L_1$ and let $\iota_1$ be the associated involution on $E$.
   So $P_2=\pi(\iota_1e)$ and $T_2$ is the image $\pi(L_2)$ of a line
   $L_2$ on $Q_2$ through $\iota_1e$.  Let $R_2$ be the ruling on $Q_2$
   containing $L_2$ and $\iota_2$ its associated involution.
   Now the tangents $T_i$, $i\ge 3$, are defined to be the images
   of the lines $L_i$, where for $k\ge 2$
   \be
      L_{2k-1} &:=& \mbox{ the line in $R_1$ through the point
                       $(\iota_2\iota_1)^{k-1}e$} \\
      L_{2k}   &:=& \mbox{ the line in $R_2$ through the point
                       $\iota_1(\iota_2\iota_1)^{k-1}e$} \ .
   \ee
   The intuition here is that the choice of the tangents $T_i$, $i\ge 3$,
   is compatible with the choices made for $T_1$ and $T_2$.

   Now suppose that the tangent sequence $(T_i)$
   {\em closes after $n$ steps}, i.e.\ $T_{2n+1}=T_1$.
   This occurs if and only if the translation $t := \iota_2\iota_1$
   is of order $n$.  Since this is independent of the starting point $P_1$,
   we have:

\begin{theorem}
   Suppose $C,C_1,C_2$ are smooth conics in a generic pencil in $\bbP_2$
   admitting a tangent sequence (in the sense above)
   which closes after $n$ steps.
   Then they admit such sequences starting with any point on $C$.
\end{theorem}

   The considerations above show slightly more:
   The existence of a closed tangent sequence is actually
   a property of the pair
   $C_1,C_2$ only, so it is not only independent of the starting point
   on $C$, but even on the choice of $C$ within the pencil generated
   by $C_1$ and $C_2$.


\section{The Poncelet condition for quadrics of revolution}
\label{section revolution}

   The aim of this section is to give an explicit formula for the
   Poncelet-$n$-condition on two quadrics of revolution.

   To begin with, let $q,q'$ be quadrics in $\bbP_1$, given by equations
   $$
      q(z,t)=az^2+bzt+ct^2, \qquad q'(z,t)=a'z^2+b'zt+c't^2 \ .
   $$
   The pair $(q,q')$ has three quasi-invariants under the
   $\PGL(2,\bbC)$-action:
   the {\it discriminants} $D:=b^2-4ac$, $D':=b'^2-4a'c'$ and the
   {\it jacobian} $J:=2(ac'+a'c)-bb'$.
   We have:

\begin{assertion}
   If $D,D'\ne 0$, then up to the $\PGL(2,\bbC)$-action
   the pair $(q,q')$ is uniquely determined
   by its invariants $D,D'$ and $J$.
\end{assertion}

\proof
   Because of $D'\ne 0$ we may assume $q'(z,t)=zt$, i.e.\
   $a'=0$, $c'=0$ and $b'=1$.  Then we have $D'=1$ and $J=-b$.
   So $b$ is determined by $J$ and the product $ac$ is determined
   by $D=b^2-4ac$.  For any $\alpha\in\bbC$ the projective transformation
   $
      (z:t)\mapsto(\alpha z: \frac t{\alpha})
   $
   leaves the invariants $D,D',J$ unchanged.  It transforms the
   quadric $q$ into
   $$
      \alpha^2 a z^2 + bzt + \frac c{\alpha^2} t^2 \ ,
   $$
   so by a suitable choice of $\alpha$ we can change $a$ to any nonzero
   value.
\qed

   In terms of the invariants $D,D',J$
   the discriminant $d(\lambda,\mu)$
   of the pencil $\lambda q+\mu q'$, $(\lambda:\mu)\in\bbP_1$ reads
   $$
      d(\lambda,\mu)=
         -\frac{\lambda^2}4 D+\frac{\lambda\mu}2 J - \frac{\mu^2}4 D' \ .
   $$

\bigskip
   Next we consider the quadrics
   $$
      Q_{a,b,c}:x^2+y^2=az^2+bzt+ct^2, \qquad a,b,c\in\bbC
   $$
   in $\bbP_3$ with homogeneous coordinates $(x:y:z:t)$.
   Such a quadric $Q_{a,b,c}$ is singular (a cone), if $b^2=4ac$.
   The intersection curve of two quadrics
   $Q_{a_1,b_1,c_1}$ and $Q_{a_2,b_2,c_2}$
   consists of two circles
   $$
      C_j:x^2+y^2=a_1z_j^2+b_1z_jt_j+c_1t_j^2
   $$
   in the planes $zt_j=tz_j$, $j=1,2$, where
   $(z_1:t_1)$ and $(z_2:t_2)$ are the roots of the equation
   $$
      (a_1-a_2)z^2+(b_1-b_2)zt+(c_1-c_2)t^2 = 0 \ .
   $$
   In case $(b_1-b_2)^2=4(a_1-a_2)(c_1-c_2)$ the two circles coincide,
   i.e.\ the quadrics touch along one circle.  Otherwise the circles
   $C_1,C_2$ meet in the two points $(1:i:0:0),(i:1:0:0)$.

\bigskip
   Now we give an explicit formula for the Poncelet-$n$-condition on two
   quadrics $Q_{a_1,b_1,c_1}$ and $Q_{a_2,b_2,c_2}$.  This is best
   expressed in terms of the invariants of the binary quadrics
   $q_{a_i,b_i,c_i}:=a_iz^2+b_izt+c_it^2$.

\begin{proposition}\label{Poncelet formula}
   Let $Q_{a_1,b_1,c_1}$ and $Q_{a_2,b_2,c_2}$ be two quadrics as above
   with reduced intersection
   curve (two circles) and let $D_1,D_2,J_{12}$ be the
   discriminants and the jacobian of the associated binary quadrics.
   Then the following statements
   are equivalent:
   \begin{paritems}{ii)}
      \item[i)] $Q_{a_1,b_1,c_1}$ and $Q_{a_2,b_2,c_2}$ are in
         Poncelet-$n$-position.
      \item[ii)] There is an integer $k$, $(k,n)=1$, such that
      $$
         \left(J_{12}(1+\cos\frac{2\pi k}n)+D_1+D_2\right)^2
         = D_1D_2\left(1-\cos\frac{2\pi k}n\right)^2
      $$
   \end{paritems}
\end{proposition}

\proof
   Let $C_1,C_2$ be the circles of intersection of the two quadrics
   and let $zt_j=tz_j$, $j=1,2$, be the associated planes.
   It will be enough to consider the case that $t_1=t_2=1$.
   We may then calculate in affine coordinates $x,y,z$.

   To begin with, let
   $$
      P_1=(x_1,y_1,z_1),\quad P_2=(x_2,y_2,z_2),\quad P_1'=(x'_1,y'_1,z_1)
   $$
   be three points such that the line $\overline{P_1P_2}$ lies on $Q_1$
   and $\overline{P_2P'_1}$ lies on $Q_2$.  Now,
   $\overline{P_1P_2}$ lies on $Q_1$ if and only if the points
   $P_1,P_2$ belong to $Q_1$, i.e.\
   $$
      x_j^2+y_j^2=a_1z_j^2+b_1z_j+c_1
   $$
   and if the point $(P_1+P_2)/2$ lies on $Q_1$, i.e.\
   $$
      \left(\frac{x_1+x_2}2\right)^2+\left(\frac{y_1+y_2}2\right)^2
      = a_1\left(\frac{z_1+z_2}2\right)^2+b_1\left(\frac{z_1+z_2}2\right)^2
      + c_1 \ .
   $$
   Because of $P_1,P_2\in Q_1$ this is equivalent to
   \begin{equation}\label{eq 1}
      x_1x_2+y_1y_2=a_1z_1z_2+\frac{b_1}2(z_1+z_2)+c_1 =: g_1 \ .
   \end{equation}
   Similarly the line $\overline{P_2P_1'}$ lies on $Q_2$ if and only if
   \begin{equation}\label{eq 2}
      x_1'x_2+y_1'y_2=a_2z_1z_2+\frac{b_2}2(z_1+z_2)+c_2 =: g_2 \ .
   \end{equation}
   Now we come to the Poncelet condition.  The automorphisms of order
   $n$ on $C_1\setminus\{(1:i:0:0),(i:1:0:0)\}\isom\bbC^*$ are
   induced by the maps
   \begin{equation}\label{eq 3}
      \rho_k:
      \left(
      \begin{array}{cc}
         x \\ y
      \end{array}
      \right)
      \mapsto
      \left(
      \begin{array}{rr}
         \cos\tfrac{2\pi k}n&-\sin\tfrac{2\pi k}n \\
         \sin\tfrac{2\pi k}n&\cos\tfrac{2\pi k}n
      \end{array}
      \right)
      \left(
      \begin{array}{cc}
         x \\ y
      \end{array}
      \right)
   \end{equation}
   where $(k,n)=1$.  So the two quadrics are in Poncelet-$n$-position
   if and only if
   $\rho_k(x_1,y_1)=(x_1',y_1')$
   for some $k$ with $(k,n)=1$.
   Without loss of generality we may assume $y_1=0$.  Then
   conditions \eqnref{eq 1}, \eqnref{eq 2} and \eqnref{eq 3} read
   \begin{eqnarray}
      x_1'&=&x_1\cos\tfrac{2\pi k}n \\
      -y_1'&=&x_1\sin\tfrac{2\pi k}n \\
      x_1 x_2&=&g_1 \label{eq x_2}   \\
      x_1'x_2+y_1'y_2&=&g_2
   \end{eqnarray}
   Inserting the first three equations into the last one we obtain
   \begin{equation}\label{eq y_2}
      g_1\cos\tfrac{2\pi k}n - x_1y_2\sin\tfrac{2\pi k}n = g_2 \ . \label{y_2}
   \end{equation}
   and
   from the equations \eqnref{eq x_2} and \eqnref{eq y_2} we get
   $$
      x_2^2+y_2^2=\left(\frac{g_1}{x_1}\right)^2+
      \left(
      \frac{g_1\cos\tfrac{2\pi k}n-g_2}{x_1\sin\tfrac{2\pi k}n}
      \right)^2 \ .
   $$
   Now, upon using the equation of $Q_1$ we can replace the left hand
   side by $a_1z_2^2+b_1z_2+c_1$ and we can substitute $x_1^2$ on the
   right hand side by $a_1z_1^2+b_1z_1+c_1$.
   In this way
   we arrive at the equation
   \begin{equation}\label{eq ponc 1}
      g_1^2-2g_1g_2\cos\tfrac{2\pi k}n+g_2^2=p\sin^2\tfrac{2\pi k}n
   \end{equation}
   where
   $p:=(a_1z_1^2+b_1z_1+c_1)(a_1z_2^2+b_1z_2+c_1)$.
   Now we can use that $z_1$ and $z_2$ are the roots of the
   equation
   $
      (a_1-a_2)z^2+(b_1-b_2)z+(c_1-c_2)=0
   $
   to write $g_1,g_2$ and $p$ in terms of the coefficients
   $a_1,b_1,c_1,a_2,b_2,c_2$.
   This finally allows to express \eqnref{eq ponc 1} in the
   invariants $D_1,D_2,J_{12}$.  We omit the details.
\qed


\section{Circles in the projective plane}

   In the previous sections we considered Poncelet properties of
   conics and quadrics related to the three-dimensional Poncelet theorem.
   Now we turn to the study of closing theorems in circle geometry.

   A {\em circle} is a conic $C\subset\bbP_2$ passing through the two
   circular points $(1:\pm i:0)$.  Its equation is of the form
   $$
      a(x^2+y^2)+2bxz+2cyz+dz^2=0 \ ,
   $$
   where $(a:b:c:d)\in\bbP_3$.  The {\em discriminant} of the circle $C$ is
   $$
      \det\left(\begin{array}{ccc} a & 0 & b \\
                                   0 & a & c \\
                                   b & c & d
                \end{array}\right)
      = a(ad-b^2-c^2) \ .
   $$
   The quadratic form
   $$
      q(C) := ad-b^2-c^2
   $$
   is the basic invariant.  A circle $C$ with $q(C)=0$ consists of two lines
   each passing through a circular point.  These are the {\em null-circles}.
   If we have $a=0$, then $C$ decomposes into two lines, one of which
   passes through the two circular points.  In this case we
   say that $C$ is a {\em line}.

   The quadratic invariant $q$ associates with each circle $C$ two basic
   surfaces in the space of circles $\bbP_3$:

\subparagraph*{The polar plane.}
   The polar plane of a circle $C$ with respect to
   the quadric $q$ is
   $$
      q(C,C') = \tfrac12(ad'+da')-(bb'+cc') \ .
   $$

\begin{proposition}[\cite{Ped57}, p.\ 32]
   We have $q(C,C')=0$ if and only if the circles $C,C'$ intersect
   orthogonally, i.e.\ if they have a point of intersection $P$ such that
   the tangent vectors of $C$ and $C'$ in $P$ are orthogonal with respect to
   the canonical bilinear form on $\bbC^2$.
\end{proposition}

\proof
   We may use affine coordinates $x,y$.
   The tangent vector of $C$ at the point $(x,y)$ is
   $$
      \vect{\d Cx \\ \d Cy} = \vect{2ax+2b \\ 2ay+2c} \ ,
   $$
   so $C$ and $C'$ intersect orthogonally if and only if there is a point
   $(x,y)$ such that the following equations are satisfied.
   \be
      a(x^2+y^2)+2bx+2cy+d &=& 0 \\
      a'(x^2+y^2)+2b'x+2c'y+d' &=& 0 \\
      (2ax+2b)(2a'x+2b')+(2ay+2c)(2a'y+2c') &=& 0
   \ee
   Elimination of $x,y$ from these equations yields $q(C,C')=0$, as claimed.
\qed

\subparagraph*{The tangent cone.}
   This is the variety of circles touching a given circle.  We have:

\begin{proposition}[\cite{Ped57}, p.\ 32] \label{tangent cone}
   a) The circles $C'$ touching a given circle $C$ are parametrized by the
   quadric surface
   $$
      Q_{C}: q(C)q(C')-q(C,C')^2=0
   $$
   in the space of circles.

   b) If $C$ is smooth, then $Q_{C}$ is a cone with top $C\in\bbP_3$.  If
   $C$ is a null-circle, then $Q_{C}$ is a double plane.
\end{proposition}

\proof
   a)
   The circles $C,C'$ touch if there is a point $(x,y)$ such that
   \be
      C(x,y)=C'(x,y)=0 \and
      dC(x,y)=\lambda\cdot dC'(x,y)
   \ee
   for some $\lambda\in\bbC^*$.  By eliminating $x,y$ and $\lambda$ from
   these equations we find
   \be
   && -\tfrac14 a^2 d'^2+\tfrac12 ada'd'-adb'^2-adc'^2+abb'd'+acc'd'-
   \tfrac14 d^2a'^2+bda'b'+cda'c'\\
   && -b^2a'd'+b^2c'^2-2bcb'c'-
   c^2a'd'+c^2b'^2 = 0
   \ ,
   \ee
   which is equivalent to
   $q(C)q(C')-q(C,C')^2=0$.

   b)
   An obvious calculation shows that $C=(a:b:c:d)$
   is a singular point of $Q_C$.
   Further, we find $\rank Q_C=1$ if $q(C)=0$.
\qed

   Now we consider the variety of all circles $C$,
   which touch two given (smooth) circles $C_1,C_2$.  This is the
   intersection curve of the two cones
   \be
      Q_{C_1}: q(C_1)q(C)-q(C_1,C)^2&=&0 \\
      Q_{C_2}: q(C_2)q(C)-q(C_2,C)^2&=&0 \ .
   \ee

\begin{proposition}\label{planes}
   a) The intersection of the cones $Q_{C_1},Q_{C_2}$ consists of two
   conics lying in the planes
   $$
      \Pi_{C_1,C_2}^{\pm}:=\sqrt{q(C_2)}q(C_1,C)\pm\sqrt{q(C_1)}q(C_2,C) \ .
   $$

   b) The intersection points of the two conics are the null-circles
   associated to the intersection points of $C_1$ and $C_2$.
\end{proposition}

   In the sequel the two conics above will be referred to as the two
   {\em families of circles touching $C_1$ and $C_2$}.

\proof
   a) The cones $Q_{C_1},Q_{C_2}$ span the pencil of quadrics
   $$
      (\lambda q(C_1)+\mu q(C_2))q(C)-\lambda q(C_1,C)^2-\mu q(C_2,C)^2
      \qquad ,(\lambda:\mu)\in\bbP_1 \ .
   $$
   This pencil contains the quadric
   $$
      q(C_2)Q_{C_1}-q(C_1)Q_{C_2}=q(C_2)q(C_1,C)^2-q(C_1)q(C_2,C)^2 \ ,
   $$
   which splits in the two planes $\Pi_{C_1,C_2}^{\pm}$.

   b) Let $C$ be a null-circle associated to $C_1,C_2$.  It touches both
   $C_1$ and $C_2$, so
   $$
      q(C_1)q(C)-q(C_1,C)^2=
      q(C_2)q(C)-q(C_2,C)^2=0 \ .
   $$
   Because of $q(C)=0$ this implies
   $q(C_1,C)=q(C_2,C)=0$, so $C$ lies on both planes
   $\Pi_{C_1,C_2}^{\pm}$.
\qed

   We will also need the following degenerate case.

\begin{lemma}\label{touching}
   If two circles $C_1,C_2$ touch each other, then the tangent cones
   $Q_{C_1},Q_{C_2}$ touch along a line.
\end{lemma}

\proof
   We may assume that the circles have equations
   $$
      C_1: x^2+y^2+z^2=0 \and C_2:a(x^2+y^2)+2bxz+2cyz+dz^2=0 \ .
   $$
   Let $P$ be the
   pencil of circles touching $C_1$ in its point of contact with $C_2$.
   Thus $P$ is just the line joining the vertices
   $p_1=(1:0:0:1)$ and $p_2=(a:b:c:d)$ of the two cones.
   We calculate the tangent
   planes of $Q_{C_1}$ and $Q_{C_2}$ in the points
   $\lambda_1 p_1+\lambda_2 p_2$ of $P$:
   \be
      Q_{C_1}(\lambda_1 p_1+\lambda_2 p_2)&=&\lambda_2Q_{C_1}p_2
         = \tfrac{\lambda_2}4(d-a:-4b:-4c:a-d) \\
      Q_{C_2}(\lambda_1 p_1+\lambda_2 p_2)&=&\lambda_1Q_{C_2}p_1
         = \\
      \rlap{$
        \frac{\lambda_1}4(-d^2-2b^2-2c^2+ad:2b(a+d):
                           2c(a+d):
                           -a^2-2b^2-2c^2+ad)
      $}\phantom{Q_{C_2}(\lambda_1 p_1+\lambda_2 p_2)}\hskip-0.5cm
   \ee
   Now $C_1,C_2$ touch each other, i.e.\ $C_2\in Q_{C_1}$, i.e.\
   $a^2+d^2+4b^2+4c^2-2ad=0$.  This shows that
   the tangent planes coincide along $P$.
\qed

   In Sect.\ \ref{section Money-Coutts}
   we will have to work with the tangent cones of three given
   circles $C_1,C_2,C_3$.  The following property of their mutual
   intersection will turn out to be crucial.

\begin{lemma}\label{common line}
   If the
   signs $\eps_{12},\eps_{13},\eps_{23}\in\{\pm1\}$ are chosen
   such that
   $\eps_{12}\eps_{13}\eps_{23}=-1$, then the planes
   $
   \Pi_{C_1,C_2}^{\eps_{12}},
   \Pi_{C_1,C_3}^{\eps_{13}},
   \Pi_{C_2,C_3}^{\eps_{23}}
   $
   have a line in common.
\end{lemma}

   Notice: Changing the sign of a square root $\sqrt{q(C_i)}$
   does not influence the condition
   $\epsilon_{12} \epsilon_{13} \epsilon_{23} = -1$.

\proof
   The pencil spanned by
   $\Pi_{C_1,C_2}^{\eps_{12}}$ and
   $\Pi_{C_1,C_3}^{\eps_{13}}$ contains the plane
   $$
      {\frac{\sqrt{q(C_2)}}{\sqrt{q(C_1)}}}\Pi_{C_1,C_3}^{\eps_{13}} -
      {\frac{\sqrt{q(C_3)}}{\sqrt{q(C_1)}}}\Pi_{C_1,C_2}^{\eps_{12}}
      =
      \eps_{13}\sqrt{q(C_2)}q(C_3,C) -
      \eps_{12}\sqrt{q(C_3)}q(C_2,C) \ ,
   $$
   which equals
   $\Pi_{C_2,C_3}^{\eps_{23}}$, if
   $\eps_{12}\eps_{13}\eps_{23}=-1$.
\qed


\section{Emch's and Steiner's theorems on circular series}
\label{section Emch Steiner}

   In this section we consider two classical theorems on circular series.
   Our aim here is to show
   that one obtains short proofs by considering an
   elliptic resp.\ rational curve underlying the closing mechanism.

   Let $C,C_1,C_2\subset\bbP_2$ be smooth circles
   in general position and let $F$ be one of the two families of circles
   touching $C_1$ and $C_2$.

   We start with a smooth circle $S_1\in F$ and choose
   one of its points of intersection with
   $C$, $P_1$ say.
   The pair $(S_1,P_1)\in F\times C$ determines a second pair as
   follows.  There are two circles in $F$ through $P_1$, namely the points
   of intersection of $F$ with the plane $\Pi_{P_1}$
   in the space of circles consisting
   of all circles through
   $P_1$.  Let $F\cap\Pi_{P_1}=\{S_1,S_2\}$.
   The second circle $S_2$ in turn
   intersects $C$ in $P_1$ and in a second point $P_2$.  We set
   $t(S_1,P_1) := (S_2,P_2)$.

\begin{theorem}[Emch \cite{Emc01}]
   Suppose $t^n(S_1,P_1)=(S_1,P_1)$ holds for some pair $(S_1,P_1)$
   such that $P_1\in C$.
   Then this holds for any such pair.
\end{theorem}

\proof
   Consider the incidence curve
   $$
      E := \{(S,P)\in F\times C\with P\in S\}
   $$
   Denoting by $\pi_1,\pi_2$ the projections onto $F$ resp.\ $C$, we
   have for $S\in F$
   $$
      \pi_2\pi_1\inverse(S) = S\cap C=Q+Q'+P+P' \ ,
   $$
   where $Q,Q'$ are the infinitely far points on $C$ and $P,P'\in C$.
   For $P\in C$ we have
   $\pi_1\pi_2\inverse(P)=F\cap \Pi_P$,
   where $\Pi_P$ is the plane in the space of circles consisting
   of the circles through $P$.
   Thus the curve $V$ is of bidegree $(4,2)$ and it decomposes as
   $$
      V=F\times \{Q\} + F\times\{Q'\} + E \ ,
   $$
   where $E$ is a curve of bidegree $(2,2)$.
   The branch points of the restriction $\pi_1\restr E$ are the circles
   $S\in F$ touching $C$.  Since $C,C_1,C_2$ are in general position,
   there are exactly four of them.
   So $E$ is smooth elliptic by Lemma \ref{smooth elliptic} below.
   The map $t$ is just the composition of the covering involutions of
   $\pi_1\restr E$ and $\pi_2\restr E$, so it is a translation on $E$.
   The assumption implies that it is of order $n$.
\qed

\begin{lemma}\label{smooth elliptic}
   Let $E\subset\bbP_1\times\bbP_1$ a curve of bidegree $(2,2)$.  Assume
   that the projections $\pi_1,\pi_2:E\to\bbP_2$ are finite and that one of
   them has at least four branch points.  Then $E$ is smooth.
\end{lemma}

\proof
   Suppose to the contrary that $E$ is singular.  First assume that $E$ is
   reducible.  Since the projections are finite, $E$ is then a sum of two
   curves $E_1,E_2$ of bidegree $(1,1)$.  By the adjunction formula
   $E_1\isom E_2\isom\bbP_1$.  But then $\pi_1$ and $\pi_2$ have only two
   branch points, namely the intersection points of $E_1$ and $E_2$, a
   contradiction.

   If $E$ is irreducible, then it has exactly one singularity and the
   normalization $\tilde E$ is a smooth rational curve.  But then by
   the Hurwitz formula the map $\tilde E\to E\to\bbP_1$ has only two branch
   points, a contradiction again.
\qed

   Next we turn to Steiner's theorem.  In contrast to the situation
   of Emch's theorem here the relevant curve for the closing process
   is rational.

   Let $C_1,C_2\subset\bbP_2$ be smooth circles in general position
   and let $F$ be a family of circles touching $C_1$ and $C_2$.

\begin{theorem}[Steiner]
   Suppose that
   there is a closed sequence of distinct
   circles $S_1,S_2,\dots,S_n,S_{n+1}=S_1$
   in $F$ such that $S_i$ touches $S_{i+1}$ for $1\le i\le n$.
   Then there are such sequences starting with any circle in $F$.
\end{theorem}

   Steiner's theorem is a consequence of the following

\begin{proposition}
   a) $V := \{(S_1,S_2)\in F\times F\with S_1\in Q_{S_2} \}$
   is a curve of bidegree $(4,4)$ and
   $$
      V = 2\cdot\Delta + V_1 + r^*V_1 \ ,
   $$
   where $\Delta\subset F\times F$ is the diagonal,
   $V_1$ is a curve of bidegree $(1,1)$ and
   $r:\bbP_1\times\bbP_1\to\bbP_1\times\bbP_1, (S_1,S_2)\mapsto(S_2,S_1)$.

   b) We have $V_1\setminus\Delta\isom r^*V_1\setminus\Delta
   \isom\bbC^*$.  Let
   \be
      t:V_1\setminus\Delta &\to& V_1\setminus\Delta
   \ee
   be defined by $t(S_1,S_2) := (S_2,S_3)$, where
   $\pi_2^{-1}(S_2)=2(S_2,S_2)+(S_1,S_2)+(S_3,S_2)$.
   Then $t$ is the multiplication by a nonzero constant.
\end{proposition}

\proof
   We may consider $V$ as a subvariety of $\bbP_3\times\bbP_3$, namely as the
   intersection of the varieties $F\times\bbP_3$, $\bbP_3\times F$ and the
   hypersurface $\{(S_1,S_2)\with S_1\in Q_{S_2}\}$ of bidegree $(2,2)$.
   Denoting by $\pi_1,\pi_2:V \to F$ the projections, we have for $S_1\in F$
   $$
      \pi_2\pi_1^{-1}(S_1) = Q_{S_1}\cap Q_{C_1} \restr{\Pi(F)} \ ,
   $$
   where $\Pi(F)\subset\bbP_3$ is the plane containing $F$.
   According to Lemma \ref{touching} the cones $Q_{S_1},Q_{C_1}$ touch
   along a line, hence the restriction to $F$ of their intersection is
   of the form
   $$
      Q_{S_1}\cap Q_{C_1} \restr{\Pi(F)} =2 S_1 + S_2 + S_2' \ ,
   $$
   where $S_2,S_2'\in F$.  Therefore $V$ contains the diagonal $\Delta$
   as a component of multiplicity two, the residual curve
   $V'$ being of bidegree $(2,2)$.

   Now let $S_0$ be one of the null-circles associated to $C_1$ and $C_2$,
   i.e.\ $S_0$ consists of the two lines joining an intersection point $P$ of
   $C_1,C_2$ with the infinitely far circle points.  The quadric $Q_{S_0}$ is
   then of rank 1, i.e.\ a double plane.  It consists of the circles
   through $P$.  The set $\pi_2\pi_1^{-1}(S_0)$ contains the circles
   through $P$ touching $C_1$ and $C_2$, so $\pi_1^{-1}(S_0)=\{(S_0,S_0)\}$.
   We conclude that the two null-circles are branch points of both
   projections $\pi_1,\pi_2$.  Because of $p_a(V')=1$ this implies that
   $V'$ is reducible.  It consists of two curves $V_1,V_2$ of bidegree
   $(1,1)$ meeting in the points $(S_0,S_0),(S_0',S_0')$, where $S_0,S_0'$
   are the two null-circles associated to $C_1$ and $C_2$.
   Because of $r^*V=V$, we have $V_2=r^*V_1$.
\qed


\section{The zig-zag theorem}
\label{section zig-zag}

   In this section we give a proof of (a complex-projective version of) the
   zig-zag theorem \cite{BlaHowHow74}
   based on the consideration of an elliptic curve in the product of two
   circles.

   For a point $P_0=(x_0:y_0:z_0:t_0)$ in $\bbP_3$ and a complex number $r$
   the quadric
   $$
      S_{P_0,r}: (xt_0-tx_0)^2+(yt_0-ty_0)^2+(zt_0-tz_0)^2=r^2t^2t_0^2
   $$
   is called the {\em sphere with center $P_0$ and radius $r$}.
   Its intersection with the infinitely far plane $t=0$ is the conic
   $x^2+y^2+z^2=0$.  A {\em circle} in $\bbP_3$ is a conic whose
   infinitely far points lie on this conic.

\bigskip
   Now let $C_1,C_2\subset\bbP_3$ be smooth
   circles and let a radius $r\in\bbC$ be fixed.
   A {\em zig-zag} for $C_1,C_2$ is a sequence of points $P_1,P_2,\dots$ such
   that for $i\ge 1$ the point $P_i$ lies on $C_1$ resp.\ $C_2$, if $i$ is
   odd resp.\ even, and such that
   $P_{i+1}\in S_{P_i,r}$.  Consecutive points $P_i,P_{i+1}$
   in a zig-zag are thought of as having constant ''distance'' $r$.
   The zig-zag is said to {\em close after $n$ steps}, if $P_{2n+1}=P_1$.

\begin{theorem}[cf.\ \cite{BlaHowHow74}]
   Let $C_1,C_2\subset\bbP_3$ be a general pair of circles and let
   $r\in\bbC$.
   If the pair $C_1,C_2$ admits a zig-zag of distance $r$
   which closes after $n$ steps, then
   it admits such zig-zags starting with any point on $C_1$.
\end{theorem}

\proof
   1) We consider the curve
   $$
      V := \{(P_1,P_2)\in C_1\times C_2 \with P_2\in S_{P_1,r}\} \ .
   $$
   It is the restriction of the hypersurface
   $\{P_2\in S_{P_1,r}\}\subset\bbP_3\times\bbP_3$ of bidegree $(2,2)$,
   so it is of bidegree $(4,4)$ in $C_1\times C_2$.
   We denote by $\pi_i:V\to C_i$, $i=1,2$, the projections and by $Q_1,Q_1'$
   resp.\ $Q_2,Q_2'$ the infinitely far points on $C_1$ resp.\ $C_2$.
   For a point $P_1\in C_1$, different from $Q_1$ and $Q_1'$, we have
   $$
      \pi_2\pi_1\inverse(P_1)=S_{P_1,r}\cap C_2=Q_2+Q_2'+P_2+P_2' \ ,
   $$
   where $P_2,P_2'\in C_2$.  Therefore $V$ contains the lines
   $$
      C_1\times\{Q_2\}, \ C_1\times\{Q_2'\}
   $$
   and by the same reasoning also the lines
   $$
      \{Q_1\}\times C_2, \ \{Q_1'\}\times C_2 \ .
   $$
   The residual curve $E$ is of bidegree $(2,2)$.

   2) Next we determine the branch points of the projection
   $\pi_1\restr E:E\to C_1$ to show that $E$ is smooth elliptic.
   Let us first consider the branch points of $\pi_1:V\to C_1$.
   These are the points $P_1\in C_1$ such that the sphere $S_{P_1,r}$
   touches $C_2$.  We may assume $C_2$ to lie in the plane $z=0$
   having equation $x^2+y^2+t^2=0$,
   so $S_{P_1,r}$ intersects the plane of $C_2$
   in the circle $S:=S_{P_1,r}(x,y,0,t)$.  The sphere $S_{P_1,r}$
   touches $C_2$ if and only if $S$ lies on the tangent cone of $C_2$.
   Evaluating this condition by means of \eqnref{tangent cone} we find
\begin{equation}\label{branch point quartic}
      -\tfrac14(x_1^2+y_1^2+z_1^2-r^2t_1^2)^2
      +\tfrac12 t_1^2(x_1^2+y_1^2+z_1^2-r^2t_1^2)
      -x_1^2t_1^2-y_1^2t_1^2-\tfrac14 t_1^4 = 0 \ .
\end{equation}
   So the branch points of $\pi_1$
   are just the intersection points of the circle
   $C_1$ and the quartic surface $Y\subset\bbP_3$ defined by
   \eqnref{branch point quartic}.
   One may think of $Y$ as the set of all points having ''distance''
   $r$ from $C_2$.
   Since $Y$ intersects the plane $t_1=0$ in the double conic
   $(x_1^2+y_1^2+z_1^2)^2=0$, the points $Q_1,Q_1'$ are both contained
   in $Y\cap C_1$ with multiplicity two.  The four remaining points
   of intersection are the branch points of $\pi_1\restr E:E\to C_1$.
   For general $C_1$ these points are distinct, so $E$ is in fact
   smooth elliptic.

   3)
   Having identified the underlying elliptic curve,
   the proof now ends in the usual way:
   We denote by
   $\iota_1,\iota_2$ the covering involutions of
   $\pi_1\restr E$ resp.\ $\pi_2\restr E$
   and by $t := \iota_2\iota_1$ the associated translation on $E$.
   The assumption that there is a zig-zag
   which closes after $n$ steps implies that $t$ is of order $n$.
   This proves the theorem.
\qed


\section{The Money-Coutts theorem}
\label{section Money-Coutts}

   Let $C_1,C_2,C_3\subset\bbP_2$ be smooth circles in general position.
   For each of the pairs $C_1,C_2$ and $C_2,C_3$
   we choose one of the two families of circles touching both circles,
   say $F_{1},F_{2}$ respectively.
   We are interested in the pairs of circles $S_1\in F_1,S_2\in F_2$
   such that $S_1$ touches $S_2$.
   To begin with, we show:

\begin{proposition}\label{elliptic curve}
   The curve
   $$
      V :=\{(S_1,S_2)\in F_{1}\times F_{2} \with S_1\in Q_{S_2} \}
   $$
   is of bidegree $(4,4)$ and we have
   $V = 2\tilde\Delta+E$,
   where $\tilde\Delta\subset F_1\times F_2$ is a curve of bidegree $(1,1)$
   and $E$ is a smooth elliptic curve.
\end{proposition}

\proof
   We may consider $V$ as the subvariety
   $$
      (F_1\times\bbP_3)\cap(\bbP_3\times F_2)
      \cap\{(S_1,S_2)\with S_1\in Q_{S_2} \}
   $$
   of $\bbP_3\times\bbP_3$.
   So for $S_1\in F_1$ we have
   $$
      \pi_2\pi_1^{-1}(S_1)=Q_{S_1}\cap Q_{C_2} \restr{\Pi(F_2)} \ ,
   $$
   where $\Pi(F_2)\subset\bbP_3$ is the plane containing $F_2$ and
   $\pi_i:V\to F_i$, $i=1,2$, are the projections.  By Lemma \ref{touching}
   the cones $Q_{S_1},Q_{C_1}$ touch along a line, so their intersection
   restricts to $\Pi(F_2)$ as
   $$
      Q_{S_1}\cap Q_{C_2} \restr{\Pi(F_2)}=2S_1'+S_2+S_2' \ ,
   $$
   where $S_1'\in F_2$ lies on the line of contact and $S_2,S_2'\in F_2$ lie
   on the residual conic.  This shows that any line
   $\{S_1\}\times F_2\subset F_1\times F_2$
   touches $V$, hence $V$ must contain a curve $\tilde\Delta$ of bidegree
   $(1,1)$ as a component of multiplicity two.  Next we consider the
   residual curve $E$, which is of bidegree $(2,2)$ and we determine the
   branch points of the projection $\pi_1\restr E:E\to F_1$.

   $\bullet$ First let $S_1$ be one of the two null-circles associated to
   $C_1$ and $C_2$.  In this case the tangent quadric $Q_{S_1}$ is a double
   plane.
   Therefore
   $
      \pi_2\pi_1^{-1}(S_1)=Q_{S_1}\cap Q_{C_2} \restr{\Pi(F_2)}
   $
   consists of two points only, so $S_1$ is a branch point of
   $\pi_1\restr E$.

   $\bullet$ Next let $S_1$ be an Apollonius circle of $C_1,C_2,C_3$, i.e.\
   a circle touching all three of them.  Further suppose that $S_1$ lies on
   $F_1$, but not on $F_2$.  Since the tangent cones
   $Q_{C_1},Q_{C_2},Q_{C_3}$ are in general position, there are exactly two
   such circles.  Since $S_1$ touches both $C_2$ and $C_3$, we have
   \be
      Q_{S_1}\cap Q_{C_2}&=&2L_2+D_2 \\
      Q_{S_1}\cap Q_{C_3}&=&2L_3+D_3 \ ,
   \ee
   where $L_2,L_3\subset\bbP_3$ are lines and $D_2,D_3\subset\bbP_3$
   are conics.  Thus set-theoretically
   $$
      Q_{S_1}\cap Q_{C_2}\cap Q_{C_3}=
         L_2\cap L_3+D_2\cap D_3+L_2\cap D_3+L_3\cap D_2
   $$
   consists of five points at most.  We claim that $S_1$ is a branch point
   of $\pi_1\restr E$.  To this end we show that
   only two of these points belong to $F_2$.
   We use the notation and the statement of
   Lemma \ref{common line}:  The conics $D_2$ and $D_3$ lie in planes
   $\Pi_{S_1,C_2}^{\eps_{12}}$ resp.\ $\Pi_{S_1,C_3}^{\eps_{13}}$, where
   $\eps_{12},\eps_{13}\in\{\pm1\}$, so $D_2\cap
   D_3\subset\Pi_{C_2,C_3}^{\eps_{23}}$
   for $\eps_{23}:=-\eps_{12}\eps_{13}$.  Then $L_2,L_3$ lie in
   $\Pi_{S_1,C_2}^{-\eps_{12}}$ resp.\ $\Pi_{S_1,C_3}^{-\eps_{13}}$, so
   $L_2\cap L_3\subset\Pi_{C_2,C_3}^{\eps_{23}}$ as well.
   Since $S_1=L_2\cap L_3$ does not belong to $F_2$ by assumption, we find
   \be
      Q_{S_1}\cap Q_{C_2}\restr{\Pi(F_2)}  &=& L_2\cap D_3 + L_3\cap D_2 \\
      Q_{S_1}\cap Q_{C_3}\restr{\Pi(F_2')} &=& L_2\cap L_3 + D_2\cap D_3 \ ,
   \ee
   where $F_2'$ is the second family of circles touching $C_1$ and $C_2$.
   This shows that $\pi_2\pi_1^{-1}(S_1)$ consists of the two points
   $L_2\cap D_3,L_3\cap D_2$.

   Applying Lemma \ref{smooth elliptic} we conclude that $E$ is a
   smooth elliptic curve.
\qed

   As before, let $F_1,F_2$ be families of circles touching $C_1,C_2$ resp.\
   $C_2,C_3$, and suppose $F_3$ is a family of circles touching $C_3$ and
   $C_1$.  By Proposition \ref{elliptic curve} we have three elliptic curves
   \be
      E_1 &\subset& F_1\times F_2 \\
      E_2 &\subset& F_2\times F_3 \\
      E_3 &\subset& F_3\times F_1 \ .
   \ee
   Let $\pi_1,\pi_1'$ resp.\ $\pi_2,\pi_2'$ resp.\ $\pi_3,\pi_3'$ denote the
   projections onto the factors.

   The next point we want to make is:

\begin{lemma}\label{same branch points}
   Given $F_1$ and $F_2$, the family $F_3$ can be chosen in such a way that
   the projections
   $\pi_1',\pi_2$ as well as
   $\pi_2',\pi_3$ and $\pi_3',\pi_1$
   have the same branch points in $F_2$ resp.\ $F_3$ resp.\ $F_1$.
\end{lemma}

\proof
   According to Lemma \ref{common line} the family $F_3$ can be chosen such
   that the planes $\Pi(F_1),\Pi(F_2),\Pi(F_3)$ have a line in common.
   The branch points of $\pi_1'$ are the two null-circles of the pair
   $C_1,C_2$ and the two Apollonius circles in $F_2$, which do not belong
   to $F_1$.  Because of our choice of $F_3$ these Apollonius circles do not
   belong to $F_3$ either.  This shows that the branch points of $\pi_1'$
   coincide with those of $\pi_2$.  The statement on the pairs
   $\pi_2',\pi_3$ and $\pi_3',\pi_1$
   follows in the same way.
\qed

   Now suppose that families $F_1,F_2,F_3$ have been chosen as above.
   The Money-Coutts theorem on the circles $C_1,C_2,C_3$ can then be stated
   as follows.

\begin{theorem}[Tyrrel-Powell \cite{TyrPow71}]
   Suppose that circles
   $S_1\in F_1$,
   $S_2\in F_2$,
   $S_3\in F_3$
   and
   $S_4\in F_1$
   are chosen such that
   $(S_1,S_2)\in E_1$,
   $(S_2,S_3)\in E_2$ and
   $(S_3,S_4)\in E_3$.
   Then it is possible to choose circles
   $S_5\in F_2$,
   $S_6\in F_3$
   and
   $S_7\in F_1$
   in such a way that again
   $(S_4,S_5)\in E_1$,
   $(S_5,S_6)\in E_2$ and
   $(S_6,S_7)\in E_3$
   and such that
   $S_7$ coincides with $S_1$.
\end{theorem}

\proof
   According to
   Lemma \ref{same branch points} the projections
   $\pi_1',\pi_2$ resp.\
   $\pi_2',\pi_3$ resp.\
   $\pi_3',\pi_1$
   have the same branch points, hence there are isomorphisms
   $\phi_1:E_1\to E_2$,
   $\phi_2:E_2\to E_3$
   and
   $\phi_3:E_3\to E_1$
   such that
   $\pi_2\circ\phi_1=\pi_1'$,
   $\pi_3\circ\phi_2=\pi_2'$
   and
   $\pi_1\circ\phi_3=\pi_3'$.
   We identify
   $E_1=E_2=E_3=:E$ and $\pi_1'=\pi_2,\pi_2'=\pi_3$ by means of $\phi_1$ and
   $\phi_2$.
   In this way $\phi_3$ is identified with an
   automorphism $t$ of $E$ such that $\pi_1\circ t=\pi_3'$.
   Since the elliptic curve $E$ is determined by the intersection
   points of $C_1,C_2$ and two Apollonius circles of $C_1,C_2,C_3$,
   for general $C_1,C_2,C_3$ the curve $E$ is general
   as well.  Therefore the automorphism $t$ is either a translation
   or an involution.
   By a suitable choice of $\phi_1$ and $\phi_2$ we can achieve that $t$ is
   actually a translation.
   We denote by $\iota_1,\iota_2,\iota_3$ the covering involutions of
   $\pi_1,\pi_2,\pi_3$ on $E$.

   We have $S_1=\pi_1(e)$ for some $e\in E$.
   Then $S_2\in\pi_1'\pi_1\inverse(S_1)$, which means
   $S_2=\pi_1'(\alpha_1e)$ where $\alpha_1\in\{1,\iota_1\}$.
   In the same way we proceed with $S_3$ and $S_4$ to get
   $$
      S_4=\pi_3'(\alpha_3\alpha_2\alpha_1e) =
      \pi_1(t\alpha_3\alpha_2\alpha_1e) \ ,
   $$
   where $\alpha_i\in\{1,\iota_i\}$ for $1\le i\le 3$.
   Now we determine the circles $S_5,S_6,S_7$ such that
   $(S_4,S_5),(S_5,S_6),(S_6,S_7)\in E$.  We find
   $$
      S_7=\pi_3'(\beta_3\beta_2\beta_1t\alpha_3\alpha_2\alpha_1e) =
      \pi_1(t\beta_3\beta_2\beta_1t\alpha_3\alpha_2\alpha_1e) \ ,
   $$
   where $\beta_i\in\{1,\iota_i\}$ for $1\le i\le 3$.
   Now we make our choice for $S_5,S_6,S_7$ resp.\ for
   $\beta_1,\beta_2,\beta_3$ in the following way.
   We let $\beta_2:=\alpha_2$, $\beta_3:=\alpha_3$ and we
   choose $\beta_1$ such that the number of subscripts $i$, $1\le i\le 3$,
   such that $\beta_i=\iota_i$ is odd.
   With these choices we obtain
   $$
      t\beta_3\beta_2\beta_1t\alpha_3\alpha_2\alpha_1e =
      \beta_3\beta_2\beta_1\alpha_3\alpha_2\alpha_1e =
      (\alpha_3\alpha_2\beta_1)^2\beta_1\alpha_1e \ .
   $$
   Since the composition of an odd number of involutions is again an
   involution, the latter expression equals $\beta_1\alpha_1e$,
   so $S_7=\pi_1(\beta_1\alpha_1e)=\pi_1(e)=S_1$.
\qed

   Finally, we aim at a Poncelet-type statement
   in the situation of the Money-Coutts theorem.  So suppose that
   $C_1,C_2,C_3$ are three circles as above and let families $F_1,F_2,F_3$
   be chosen as before.  As in the proof of the Money-Coutts theorem
   the elliptic curve defining the contact relation
   between circles of two families will be denoted by $E$.

   Let $(S_i)_{i\ge 1}$ be a sequence of circles such that
   $(S_i,S_{i+1})\in E$ for $i\ge 1$ and
   $$
      S_{3k+l}\in F_l \qquad\mbox{ for } k\ge 0 \mbox{ and } 1\le l\le 3 \ .
   $$
   If it were to happen that $S_{3n+1}=S_1$ for some integer $n\ge 1$, then
   the sequence is said to {\em close after $n$ steps}.

\begin{theorem}
   Suppose there are more than $2^{3n+2}$ sequences
   $(S_i)_{i\ge 1}$ closing after $n$ steps.
   Then there are infinitely many sequences
   closing after $n$ steps,
   starting with any given
   circle $S'_1\in F_1$.
\end{theorem}

\proof
   Let $S_1=\pi_1(e)$, $e\in E$.  The circles $S_i$, $i\ge 2$, are determined
   by repeatedly choosing automorphisms in $\{1,\iota_1,\iota_2,\iota_3\}$.
   So for $k\ge 1$ we have
   \begin{eqnarray*}
      S_{3k+1} &=&
       \pi_1\(\prod_{i=1}^k t\alpha_3^{(i)}\alpha_2^{(i)}\alpha_1^{(i)}\) \\
      S_{3k+2} &=&
       \pi_1'\(\alpha_1^{(k+1)}\prod_{i=1}^k
t\alpha_3^{(i)}\alpha_2^{(i)}\alpha_1^{(i)}\) \\
      S_{3k+3} &=&
       \pi_2'\(\alpha_2^{(k+1)}\alpha_1^{(k+1)}\prod_{i=1}^k
t\alpha_3^{(i)}\alpha_2^{(i)}\alpha_1^{(i)}\) \ ,
   \end{eqnarray*}
   where $t:E\to E$ is a translation and
   $\alpha_j^{(i)}\in\{1,\iota_j\}$ for $1\le j\le 3$ and $i\ge 1$.
   By assumption the sequence $(S_i)$ closes  after $n$ steps, i.e.\
   $S_{3n+1}=S_1$, i.e.\ $\pi_1(\gamma e)=\pi_1(e)$, where $\gamma$ is the
   automorphism
   $\prod_{i=1}^n t\alpha_3^{(i)}\alpha_2^{(i)}\alpha_1^{(i)}$.
   The map $\gamma$ is
   an involution or a translation on $E$,
   so replacing $e$ by $\iota_1 e$ we may assume that $\gamma$
   is an involution.
   We have $\gamma(e)=e$ or $\gamma(e)=\iota_1(e)$.
   In the first case $e$ is one of the four fixed points of $\gamma$.
   Since there are $2^{3n}$ choices for the automorphisms
   $\alpha_j^{(i)}$, $1\le j\le 3$, $1\le i\le n$, this case occurs
   for at most $2^{3n+2}$ sequences.  In the second case
   we conclude $\gamma=\iota_1$, so if we are given any circle $S_1'\in F_1$,
   $S_1'=\pi(e')$, then we define a sequence $(S_1')_{i\ge 1}$ by means
   of the automorphisms $\alpha_j^{(i)}$ and find
   $$
      S_{3n+1}=\pi_1(\gamma(e'))=\pi_1(\iota_1(e'))=\pi_1(e')=S_1 \ .
   $$
\qed


\bigskip\bigskip\bigskip
\noindent
   W.\ Barth, Th.\ Bauer \\
   Mathematisches Institut der Universit\"at\\
   Bismarckstra\ss e 1$\frac12$\\
   D-91054 Erlangen


\end{document}